\documentclass[12pt,twoside]{article}
\usepackage{times}
\usepackage{amsmath,amssymb}
\usepackage{amsthm}
\usepackage[pagewise]{lineno}
\usepackage{mathrsfs}
\usepackage{color}
\usepackage{exscale}
\usepackage{relsize}
\usepackage[title]{appendix}
\usepackage{color}
\usepackage{fancyhdr}
\allowdisplaybreaks

\newcommand{\be}{\begin{equation}}
\newcommand{\ee}{\end{equation}}
\newcommand{\bes}{\begin{equation*}}
\newcommand{\ees}{\end{equation*}}
\newcommand{\bali}{\begin{aligned}}
\newcommand{\eali}{\end{aligned}}

\newcommand{\bR}{{\mathbb R}}
\newcommand{\bN}{{\mathbb N}}

\def\nn{\nonumber}

\def\ve{\varepsilon}
\def\la{\lambda}

\def\t{\tilde}
\def\q{\quad}

\def\th{\theta}

\def\Dl{\Delta}
\def\ve{\varepsilon}

\def\lt{\left}
\def\les{\lesssim}
\def\rt{\right}

\def\i{\infty}

\def \ls{\lesssim}
\def\p{\partial}
\def\f{\frac}
\def\na{\nabla}
\def\al{\alpha}

\def\O{\Omega}

\def\s{\sqrt}

 \allowdisplaybreaks

\begin{document}

 \footskip=15pt
 \footnotesep=2pt
\let\oldsection\section
\renewcommand\section{\setcounter{equation}{0}\oldsection}
\renewcommand\thesection{\arabic{section}}
\renewcommand\theequation{\thesection.\arabic{equation}}
\newtheorem{theorem}{Theorem}[section]
\newtheorem{lemma}[theorem]{Lemma}
\newtheorem{proposition}[theorem]{Proposition}
\newtheorem{assumption}[theorem]{Assumption}
\newtheorem{remark}[theorem]{Remark}
\newtheorem{corollary}[theorem]{Corollary}

\title{Asymptotic properties of generalized D-solutions to the stationary axially symmetric Navier-Stokes equations}

\author{Zijin Li $^{a,}$\footnote{E-mail: zijinli@nuist.edu.cn} ,\quad Xinghong Pan$^{b,}$\footnote{Corresponding author email: xinghong{\_}87@nuaa.edu.cn}\vspace{0.5cm}\\
 \footnotesize $^a$School of Mathematics and Statistics, Nanjing University of Information Science and Technology,\\
  \footnotesize Nanjing 210044, China.\\
  \footnotesize $^b$Department of Mathematics, Nanjing University of Aeronautics and Astronautics, Nanjing 211106, China.\\
\vspace{0.5cm}
}

\date{}

\maketitle

\centerline {\bf Abstract}
In this paper, we derive asymptotic properties of both the velocity and the vorticity fields to the 3-dimensional axially symmetric Navier-Stokes equations at infinity under the generalized D-solution assumption $\int_{\mathbb{R}^3}|\nabla u|^qdx<\infty$  for $2<q<\infty$. We do not impose any zero or nonzero constant vector asymptotic assumption to the solution at infinity. Our results generalize those in \cite{CJ:2009JMFM,Ws:2018JMFM,CPZ2018} where the authors focused on the case $q=2$  and the velocity field approaches zero at infinity. Meanwhile, when $q\to 2_+$ and the velocity field approaches zero at infinity, our results coincide with the results in \cite{CJ:2009JMFM,Ws:2018JMFM,CPZ2018}.
\vskip 0.3 true cm

\vskip 0.3 true cm

\vskip 0.3 true cm


\section{Introduction}
In this paper, we consider asymptotic properties of smooth solution to the stationary 3D incompressible Navier-Stokes equations
\be\label{NS}
\lt\{
\begin{aligned}
&u\cdot\na u+\na p-\Dl u=0,\q x\in\bR^3;\\
&\na\cdot u=0,
\end{aligned}
\rt.
\ee
with generalized finite Dirichlet integral
\be\label{gdc0}
\int_{\bR^3}|\na u|^qdx<+\i, \q \text{for}\ 2<q<\i.
\ee
Here $u(x)\in\bR^3,p(x)\in\bR$ represent the velocity vector and the scalar pressure. Physically \eqref{NS}$_1$ represents the conservation of momentum while \eqref{NS}$_2$ shows the conservation of mass. We can also consider the same problem in an exterior domain $\O\in\bR^3$ with non-slip boundary conditions, where the complement of $\O$ is a compact axially symmetric domain, and all the results in the following can be extended to this case. However, for simplicity, we only deal with the whole space case in this paper.

The existence of weak solutions to \eqref{NS} is due to Leray \cite{Leray1933}, where he constructed a weak solution with the velocity prescribed to be a constant vector at infinity and zero at the boundary of an exterior domain. Also Leray's weak solution satisfies the bounded Dirichlet integral $\int_{\bR^3}|\na u|^2dx<+\i$. A weak solution satisfying the bounded Dirichlet integral is often referred to as ``D-solution''. See also \cite{Loa:1959UMN,Fh:1961JFSUT}. The smoothness of D-solutions is easy to prove by the properties of elliptic partial differential equations. However, the uniqueness of D-solutions has been a long and old open problem. See \cite{Galdi2011,Chae2014, KTW:2017JFA,Sg:2016NONL, CPZ2018-1, CW:2019CVPDE, PL:2019ARXIV} for some recent progress in this aspect.

An interesting and natural question is that whether weak solutions with generalized bounded Dirichlet integral
\be
\int_{\bR^3}|\na u|^qdx<+\i,\q q\neq 2,\nn
\ee
exists or not. If there exists a constant vector $u_{\i}$ such that
\be
\lim\limits_{|x|\rightarrow\i}u=u_{\i},\nn
\ee
this problem has already been investigated by several authors. For the case $q\in(2,\infty)$, the answer is positive and quite trivial. On the other  hand if $q\in(1,2)$, this situation seems to be more involved and in some situation it is hard to get the existence theorem. See \cite{GP:1991,KS:1993RM,KY:1998MANN} and references therein. Since the existence theorem  for the case  $q\in (1,2)$ is more complicated and incomplete, it is reasonable to assume that $q\in [2,\i)$.

We define a weak solution of \eqref{NS} with \eqref{gdc0} ($2\leq q<+\i$) by ``generalized D-solution''. In this paper, we restrict $q\in (2,+\i)$.

In 2 dimensional exterior domain $\O$, for the investigation of asymptotic properties of D-solution, Gilbarg-Weinberger \cite{Gilbarg1978} showed  if $u$ solves the 2D stationary Navier-Stokes equations with finite Dirichlet integral condition $\int_{\O}|\na u|^2dx<+\i$, then there exists a constant vector $u_\infty\in\mathbb{R}^2$ s.t.
\be
\lim_{r\to\infty}\int_0^{2\pi}|u(r,\th)-u_\infty|d\th=0,\nn
\ee
with the following decay estimate of vorticity:
\be
w(r,\th)=o(r^{-3/4})\quad\text{ uniformly in }\th\in[0,2\pi]\text{ as }r\to\infty,\nn
\ee
where $w:=\p_{x_2}u^1-\p_{x_1}u^2$. See also \cite{Ajc:1988ACTA,KPR:2019ARMA} for some related improvements. Recently Kozono-Terasawa-Wakasugi \cite{Kozono2019} showed that solutions of \eqref{NS} in 2D space with \eqref{gdc0} ($2<q<+\i$) satisfy a priori estimates $u(x)=o(|x|^{1-\f{2}{q}})$ and $w(x)=o(|x|^{-\left(\frac{1}{q}+\frac{1}{q^2}\right)})$ as $|x|\to\infty$.

Recently, research on the Liouville theorem of (generalized) D-solutions to the Navier-Stokes equations becomes a more and more popular topic and sufficiently fast decay of the solution at infinity is a guarantee of proving the Liouville-type theorem. If the domain $\O$ is $\mathbb{R}^2$, by applying the maximum principle of the 2D vorticity equation
\[
\Dl w-u\cdot\nabla w=0,
\]
any uniform decay of $w$ at infinity actually indicates that $w\equiv 0$. Then by Biot-Savart law, we have $-\Dl u=\na\times w=0$ which implies that $u\equiv c$ if $u$ is sublinear growth with respect to the distance to the origin. So, in 2-dimensional spaces, the generalized D-solution assumption \eqref{gdc0} implies the solution of \eqref{NS} is trivial.  However, in 3-dimensional spaces, due to the appearance of the vortex stretching term in the 3D vorticity equations, the vorticity does not satisfy the maximum principle any longer. Thus the related 3D Liouville-type problem remains open, even in the axially symmetric case. Nevertheless, a good a priori asymptotic estimate for the solution itself is significant and surely will be a cornerstone to solve the problem.

In this paper, we consider the asymptotic properties of axially symmetric generalized D-solutions to \eqref{NS} with \eqref{gdc0} in 3 dimensional space.

In the cylindrical coordinate $(r,\th,z)$, we have $x=(x_1,x_2,x_3)=(r\cos\th, r\sin\th, z)$ and a solution $u$ of \eqref{NS} is called axially symmetric if all the 3 directions of $u$ in the cylindrical coordinate do not depend on $\th$, i.e.
\[
u=u^r(r,z)e_r+u^\th(r,z)e_\th+u^z(r,z)e_z,
\]
where the basis vector $e_r$, $e_\th$ and $e_z$ are
\be\label{unitvec}
 e_r=\left(\frac{x_1}{r},\frac{x_2}{r},0\right),\quad e_\th=\left(-\frac{x_2}{r},\frac{x_1}{r},0\right),\quad e_z=\left(0,0,1\right).
\ee
Later on, we will simply denote $u=(u^r,u^\th, u^z)$. We can derive the stationary Navier-Stokes equations in cylindrical coordinate:
\begin{equation}\label{Main}
\left\{
\begin{aligned}
&(b\cdot\nabla)u^r -\frac{(u^\th)^2}{r}+\p_r p=\left(\Delta-\frac{1}{r^2}\right)u^r, \\
&(b\cdot\nabla) u^\th+\frac{u^\th u^r}{r}=\left(\Delta-\frac{1}{r^2}\right)u^\th , \\
&(b\cdot\nabla)u^z+\p_z p=\Delta u^z ,                                    \\
&b=u^re_r+u^ze_z,\q \nabla\cdot b=\p_ru^r+\frac{u^r}{r}+\p_zu^z=0.
\end{aligned}
\right.
\end{equation}
We also write the vorticity field $w=\nabla\times u$ in cylindrical coordinate:
\[
w=w^r(r,z)e_r+w^\th(r,z)e_\th+w^z(r,z)e_z,
\]
where
\[
w^r=-\p_zu^\th,\quad w^\th=\p_zu^r-\p_ru^z,\quad w^z=\p_ru^\th+\frac{u^\th}{r},
\]
and they satisfy
\be\label{VEQ}
\lt\{
\begin{aligned}
&(b\cdot\na)w^r-\left(\Dl-\f{1}{r^2}\right)w^r-(w^r\p_r+w^z\p_z)u^r=0,\\
&(b\cdot\na)w^\th-\left(\Dl-\f{1}{r^2}\right)w^\th-\f{u^r}{r}w^\th-\f{1}{r}\p_z(u^\th)^2=0,\\
&(b\cdot\na)w^z-\Dl w^z-(w^r\p_r+w^z\p_z)u^z=0.
\end{aligned}
\rt.
\ee
Recent years, a lot of studies have been devoted to the asymptotic behavior of 3D axially symmetric solution for \eqref{NS} with \eqref{gdc0} for $q=2$ and $u$ that approaches zero at infinity. We refer readers to \cite{CJ:2009JMFM,Ws:2018JMFM,CPZ2018,CPZ2018-1}, etc.. And to the best of our knowledge, the optimal results for the decay of $u$ and $w$ when $r\to\infty$ are
\be\label{QE2}
\begin{split}
|u(r,z)|&\lesssim\left(\frac{\log r}{r}\right)^{1/2};\\
|w^\th(r,z)|\lesssim\frac{(\log r)^{3/4}}{r^{5/4}},&\quad|w^r(r,z)|+|w^z(r,z)|\lesssim\frac{(\log r)^{11/8}}{r^{9/8}}.
\end{split}
\ee

Since we focus on the asymptotic properties of generalized D-solutions with \eqref{gdc0}, a larger $q$ implies a weaker assumption on the decay property of $\nabla u$ at far-field. In addition, we will not even generally assume $\lim_{|x|\to\infty}u(x)= u_\i$  for some zero or nonzero constant vector, since it is inappropriate when $q\geq3$, where $u$ may increase when $r$ tends to infinity. Meanwhile we will prove $u$ converges to a constant vector field as $r\to\infty$ when $2<q<3$. Our method is based on the scaling property of the NS system and the {\em Brezis-Wainger} inequality.

 We use $\mathcal{D}_{\lambda,l}$ and  $\t{\mathcal{D}}_{\lambda,l}$ to denote the following 2-dimensional domains
\bes
\mathcal{D}_{\lambda,l}:=\{(r,z):\,\,\lambda/4\leq r\leq4\lambda,\ |z-l|\leq 4\lambda\}
\ees
and
\be
\t{\mathcal{D}}_{\lambda,l}:=\{(r,z):\,\,\lambda/8\leq r\leq8\lambda,\ |z-l|\leq 8\lambda\},\nn
\ee
respectively. When $\lambda=1$ and $l=0$, we simply write $\mathcal{D}$, $\tilde{\mathcal{D}}$ instead of $\mathcal{D}_{1,0}$, $\tilde{\mathcal{D}}_{1,0}$. The following is our main result for the velocity:
\begin{theorem}\label{thm1}
Let  $u$ be a smooth axially symmetric solution to the Navier-Stokes equations satisfying \eqref{gdc0}. Then the oscillation of $u$ satisfies the following a priori bound
\be
\mathop{osc}_{(r,z)\in \mathcal{D}_{R,l}} u\leq CR^{1-3/q},\nn
\ee
where $C$ is a constant independent of $R$ and $l$. Furthermore,

(i) if $2<q<3$, there exists a constant $u_\infty^z$ such that
\be\label{uqs3}
|u(r,z)-u_\infty^ze_z|\leq Cr^{1-3/q},
\ee
where $C$ is independent of $r$ and $z$, and $e_z$ is the unit vector defined in \eqref{unitvec};

(ii) if $q=3$, $u$ satisfies the following ``log-growing" estimate: for $r>r_0>0$,
\be\label{uqe3}
|u(r,z)-u(r_0,z)|\leq\,C\log\lt(\f{r}{r_0}+e\rt),
\ee
where $C$ is independent of $r_0$, $r$ and $z$;

(iii) if $q>3$, $u$ satisfies the ``power-growing" estimate: for $r>r_0\geq0$,
\be\label{uqb3}
|u(r,z)-u(r_0,z)|\leq\,Cr^{1-3/q},
\ee
where $C$ is independent of $r_0$, $r$ and $z$.
\end{theorem}


\begin{remark}
Since $u$ is a smooth axially symmetric solution, we have $(u^r,\ u^\th)|_{r=0}=0$. See \cite{LW:2009JMA}. Therefore in the item (iii) of Theorem \ref{thm1}, \eqref{uqb3} indicates
\be
\sup_{z\in\mathbb{R}}\Big(|u^r(r,z)|+|u^\th(r,z)|+|u^z(r,z)-u^z(0,z)|\Big)\leq\,Cr^{1-3/q}\nn
\ee
if we choose $r_0=0$.
\end{remark}

Next we discuss the asymptotic properties of the vorticity when $r\to\infty$. Under the generalized D-condition
\[
\int_{\mathbb{R}^3}|\nabla u|^qdx<\infty,\quad q\geq 3,
\]
instead of proving a uniform bound of $u$, it seems that we can only give oscillation estimates in \eqref{uqe3} and \eqref{uqb3} by using Morrey embedding. Therefore, to derive an asymptotic behavior for $u$ that is uniformly with $z$, which will be applied to derive the asymptotic behavior of the vorticity, we need a supremum assumption on $u(r_0,z)$.
\begin{assumption}\label{assum}
If  $q\in (3,\infty)$, there exists $r_0\geq 0$ such that
\[
\sup\limits_{z\in\bR}|u(r_0,z)|\leq C({r_0});
\]

If $q=3$, there exists $r_0>0$ such that
\[
\sup\limits_{z\in\bR}|u(r_0,z)|\leq C({r_0}),
\]
where $C({r_0})$ is a constant, independent of $z$.
\end{assumption}

We have the following theorem concerning the asymptotic properties of the vorticity.

\begin{theorem}\label{thm2}
Let $u$ be a smooth axially symmetric solution to the Navier-Stokes equations \eqref{NS} with \eqref{gdc0} and $w=\nabla\times u$ be the related vorticity. Denote by $\al^-$ a positive constant which is  smaller than but close to $\al$. Then under Assumption \ref{assum}, we have

\noindent{\bf Case I:} $q\in(2,+\i)$ and $u_\i^z=0$.
\be\label{westimate1}
\begin{split}
|w^\th(r,z)|&= O\lt(r^{-(\f{1}{q}+\f{3}{q^2})^-}\rt),\\
|(w^r(r,z),w^z(r,z))|&=\lt\{
\begin{array}{ll}
 O\lt(r^{-(\f{1}{q}+\f{1}{q^2}+\f{3}{q^3})^-}\rt),&\text{ for }\ q\in(2,\f{13+\s{73}}{8});\\
 O\lt(r^{-(\f{15}{2q^2}-\f{1}{q})^-}\rt),&\text{ for }\ q\in[\f{13+\s{73}}{8},3);\\
 O\lt(r^{-(\f{1}{q-1})^-}\rt), &\text{ for }\ q\in[3,\infty),
\end{array}
\rt.
\end{split}
\ee
as $r\rightarrow+\i$.

\noindent{\bf Case II:} $q\in(2,3)$ and $u_\i^z\neq 0$.
\be\label{westimate2}
\begin{split}
|w^\th(r,z)|&= O\lt(r^{-(\f{2}{q})^-}\rt),\\
|(w^r(r,z),w^z(r,z))|&=O\lt(r^{-(\f{1}{2q}+\f{3}{q^2})^-}\rt),
\end{split}
\ee
as $r\rightarrow+\i$.
\end{theorem}


\begin{remark}
We mention here that when $q\to 2_+$ and $u_\i^z=0$, our results in Theorem \ref{thm1} and Theorem \ref{thm2} match estimates in \eqref{QE2}, except for some extra ``$log$s", owing to the critical Sobolev imbedding.
\end{remark}


\begin{remark}
When $q\in(2,3)$, estimates of the vorticity in \eqref{westimate2} are not as good as those in \eqref{westimate1} in which $u$ approaches zero at infinity. It seems strange since if $u$ approaches a non-zero constant vector at infinity, the linearized system of the Navier-Stokes equations is the Oseen system whose solutions have better decay rate at the far-field than those of the linear Stokes system. Indeed, under the assumption \eqref{gdc0} with $q=2$, the decay rate of solutions to \eqref{NS} in the case that $u$ approaches a non-zero constant vector at infinity will be better than the case that $u$ approaches zero. However, in the situation that $q>2$, it is hard to deduce a similar result. The reason is: in the case $q=2$, the nonlinear term can be regarded as a perturbation of the linear Oseen equation due to a multiplier theorem by Lizorkin \cite{Lizorkin1963} (see \cite{Galdi2011} for more details). It seems that $q=2$ is an admissible maximum in this method of perturbation, and any number $q>2$ will make the nonlinear term affect the linear Oseen equations extensively.
\end{remark}

This paper is organized as follows, in Section \ref{SEC2} we investigate the asymptotic properties of the velocity field and prove Theorem \ref{thm1}. Section \ref{SEC3} is devoted to the proof of Theorem \ref{thm2} which describes the asymptotic properties of the vorticity.

Throughout this paper, $C(c_1,c_2,...,c_n)$ denotes a positive constant depending on $c_1,\,c_2,\,...$ $\,c_n$ which may be different from line to line. For a domain $\Omega\subset\mathbb{R}^3$, $1\leq p\leq\infty$ and $k\in\mathbb{N}$, $L^p(\Omega)$ denotes the usual Lebesgue space with norm
\[
\|f\|_{L^p(\Omega)}:=
\lt\{
\begin{aligned}
&\left(\int_{\Omega}|f(x)|^pdx\right)^{1/p},\quad 1\leq p<\infty,\\
&\mathop{ess sup}_{x\in\Omega}|f(x)|,\quad\quad\quad\quad p=\infty,\\
\end{aligned}
\rt.
\]
while $W^{k,p}(\Omega)$ denotes the usual Sobolev space with its norm
\[
\begin{split}
\|f\|_{W^{k,p}(\Omega)}:=&\sum_{0\leq|L|\leq k}\|\nabla^L f\|_{L^p(\Omega)},\\
\end{split}
\]
and we simply use $H^k(\O)$ to denote the Sobolev space when $p=2$. 
We also apply $A\lesssim B$ to denote $A\leq CB$. Meanwhile, $A\sim B$ means both $A\lesssim B$ and $B\lesssim A$.

\section{Asymptotic behavior of $\boldsymbol{u}$: proof of Theorem \ref{thm1}}\label{SEC2}
Now we consider the oscillation of the solution to \eqref{NS} with \eqref{gdc0} in the domain $\mathcal{D}_{R,l}$.

\begin{lemma}[dyadic oscillation estimate]\label{DOE}
Suppose $u$ is a smooth solution of the axially symmetric Navier-Stokes equations with the generalized finite Dirichlet integral
\be\label{QDCond}
\int_{\mathbb{R}^3}|\nabla u(x)|^qdx<\infty,\q \text{for}\q 2<q<+\i.
\ee
Then the oscillation of $u$ in the domain $\mathcal{D}_{R,l}$ satisfies the following upper bound
\be\label{OSCEST}
\mathop{osc}_{\mathcal{D}_{R,l}} u\leq CR^{1-3/q},
\ee
where $C$ is a constant which is independent of $R$ and $l$.
\end{lemma}
\begin{proof}
We prove this lemma by using the scaling invariance of the Navier-Stokes equations and the  embedding theorem of Morrey. We consider the scaled solution
\be
\tilde{u}(\tilde{r},\tilde{z})=Ru(R\tilde{r},l+R\tilde{z})\nn
\ee
which is also an axially symmetric solution to the Navier-Stokes equations. We may regard $\tilde{u}$ as a two-variable function of the scaled variables $\tilde{r}$ and $\tilde{z}$ in the following two dimensional domain $\t{\mathcal{D}}$. By the imbedding theorem of Morrey (see e.g. the proof of \cite{GT1998}, Theorem 7.17), it follows that, for any $\,(\tilde{r}_1,\tilde{z}_1),\,(\tilde{r}_2,\tilde{z}_2)\in\mathcal{D}$,
\be\label{SCALEDEST}
|\tilde{u}(\tilde{r}_1,\tilde{z}_1)-\tilde{u}(\tilde{r}_2,\tilde{z}_2)|\leq C\left(\int_{\t{\mathcal{D}}}|\tilde{\nabla}\tilde{u}|^qd\tilde{r}d\tilde{z}\right)^{1/q},
\ee
where $\tilde{\nabla}=(\p_{\tilde{r}},\p_{\tilde{z}})$ and $C$ is a constant independent of $(\tilde{r}_1,\tilde{z}_1),\,(\tilde{r}_2,\tilde{z}_2)$. Now we can scale the inequality \eqref{SCALEDEST} back to the original solution $u$ and denote
\[
r=R\tilde{r},\quad z=l+R\tilde{z},\quad \bar{\nabla}=(\p_r,\p_z),
\]
and
\[
r_i=R\tilde{r}_i,\quad z_i=l+R\tilde{z}_i,\quad\text{for }\quad i=1,2,
\]
then we arrive that, $\forall\,\,({r}_1,{z}_1),\,({r}_2,{z}_2)\in\mathcal{D}_{R,l}$,
\bes
\begin{split}
R|u(r_1,z_1)-u(r_2,z_2)|&\leq C R^{2-2/q}\left(\int_{\t{\mathcal{D}}_{R,l}}|\nabla u|^qd rd z\right)^{1/q}\\
&\leq C R^{2-3/q}\left(\int_{\t{\mathcal{D}}_{R,l}}|\nabla u|^q rd rd z\right)^{1/q}.\\
\end{split}
\ees
By \eqref{QDCond}, one derives
\be\label{OSCEST1}
|u(r_1,z_1)-u(r_2,z_2)|\leq C R^{1-3/q}.
\ee
Finally, the estimate \eqref{OSCEST} holds by taking the supremum of the left-hand-side with respect to $(r_1,z_1),\,(r_2,z_2)\in \mathcal{D}_{R,l}$.
\end{proof}

Moreover, we have the following further considerations:

\subsection{Case $\boldsymbol{2<q<3}$}

\begin{proposition}
Under the same conditions as those in \textbf{Lemma \ref{DOE}} with $q\in(2,3)$, there exists a constant $u^z_\infty\in\mathbb{R}$ such that
\be
\max\Big\{|u^r(r,z)|,\,|u^\th(r,z)|,\,|u^z(r,z)-u^z_\infty|\Big\}=O(r^{1-3/q}),\quad\text{as }r\to\infty\nn
\ee
uniformly with $z\in\mathbb{R}$.
\end{proposition}
\begin{proof}
First we prove the following claim.\\
\noindent \textbf{Claim}:\\
There exists a constant vector
\bes
u_\infty=u^r_\infty e_r+u^\th_\infty e_\th+u^z_\infty e_z
\ees
such that
\be
|u(r,z)-u_\infty|=O(r^{1-3/q}),\quad\text{as }r\to\infty,\,\nn
\ee
uniformly with $z\in\mathbb{R}$. Here $e_r$, $e_\th$ and $e_z$ are unit vectors defined in \eqref{unitvec}.

$\forall z\in\mathbb{R}$, a vector field $u_{\infty}(z)$ is defined by
\be\label{Cor2.6lim-2}
u_{\infty}(z):=\lim_{n\to\infty}u(2^n,z).
\ee
This limit exists because for any $n_1>n_2\geq n$,
\be
\begin{split}
|u(2^{n_1},z)-u(2^{n_2},z)|&=\left|\sum_{i=n_2+1}^{n_1}(u(2^i,z)-u(2^{i-1},z))\right|\\
&\leq\sum_{i=n_2+1}^{n_1}\left|u(2^i,z)-u(2^{i-1},z)\right|\\
&\leq\sum_{i=n}^\infty C2^{i(1-3/q)}\\
&=C2^{n(1-3/q)}\to 0,\quad\text{ as }\quad n\to\infty.
\end{split}\nn
\ee
Here the third line follows from the oscillation estimate \eqref{OSCEST} where the constant $C$ is independent of $z$. So $\{u(2^n,z)\}_{n=1}^\infty$ is a Cauchy sequence which indicates that \eqref{Cor2.6lim-2} is well-defined and valid.

 Now we show that actually $u_{\infty}(z)$ is independent of $z$, therefore $u_{\infty}(z)=u_\i$ is a constant vector. The reason is: $\forall z_{1},\,z_{2}\in\mathbb{R}$ and $z_{1}\neq z_{2}$,
\be\label{2.13}
\begin{split}
&\q|u_{\infty}(z_{1})-u_{\infty}(z_{2})|\\
&\leq|u(2^n,z_{1})-u_{\infty}(z_{1})|+|u(2^n,z_{1})-u(2^n,z_{2})|+|u(2^n,z_{2})-u_{\infty}(z_{2})|.\\
\end{split}
\ee
$\forall\ve>0$,  by the definition \eqref{Cor2.6lim-2}, there exists an $n_0\in\mathbb{N}$ such that $\forall n>n_0$, it follows that
\be
\max\Big\{|u(2^n,z_{1})-u_{\infty}(z_{1})|\,,\,|u(2^n,z_{2})-u_{\infty}(z_{2})|\Big\}<\frac{\ve}{3}.\nn
\ee
Meanwhile, there exists an $n_0'\in\mathbb{N}$ such that $\forall n>n_0'$, $(2^n,z_{1})$ and $(2^n,z_{2})$ both belong to $\mathcal{D}_{2^{n},0}$. Now according to \eqref{OSCEST}, we arrive that
\bes
|u(2^n,z_{1})-u(2^n,z_{2})|\leq C2^{n(1-3/q)},
\ees
where $C$ is independent of $z_1$ and $z_2$. Therefore, by choosing
\be
n>\max\left\{n_0,\,n_0',\,\frac{\log\left(\frac{3C}{\ve}\right)}{\left(\frac{3}{q}-1\right)\log 2}\right\},\nn
\ee
\eqref{2.13} leads to
\be
|u_{\infty}(z_{1})-u_{\infty}(z_{2})|<\ve,\nn
\ee
which implies the constancy of $u_\infty(z)$ by choosing $\ve\to 0_{+}$. Below we use the constant $u_\i$ instead of $u_\i(z)$ for simplicity. Finally, for fixed $r>1$ and $z\in\mathbb{R}$, there exists an $n_1\in\mathbb{N}$ such that $2^{n_1}\leq r<2^{n_1+1}$. Consider the oscillation estimate \eqref{OSCEST} in domain  $\mathcal{D}_{r,z}$, one has
\be
|u(2^{n_1},z)-u(r,z)|\leq Cr^{1-3/q}.\nn
\ee
According to \eqref{Cor2.6lim-2}, there exists an $n_2\in\mathbb{N}$ (we assume $n_2>n_1$ without loss of generality) uniformly with respect to $z$, such that
\be
|u(2^{n_2},z)-u_\infty|<C r^{1-3/q}.\nn
\ee
Hence
\be
\begin{split}
|u(r,z)-u_\infty|\leq\,&\left|u(r,z)-u(2^{n_1},z)\right|+\left|u(2^{n_2},z)-u_\infty\right|+\sum_{i=n_1}^{n_2-1}\left|u(2^i,z)-u(2^{i+1},z)\right|\\
\leq&Cr^{1-3/q}+Cr^{1-3/q}+C 2^{n_1(1-3/q)}\\
\leq&Cr^{1-3/q}.
\end{split}\nn
\ee
 Here we have applied the oscillation estimate \eqref{OSCEST} in domain $\mathcal{D}_{2^{i},z}$ and $2^{n_1}\leq r<2^{n_1+1}$ to handle the last term above. This proves the Claim.
\vskip 0.1 cm
\noindent \textbf{Proof of } $\boldsymbol{u^r_\i=u^\th_\i=0.}$
\vskip 0.1 cm
Finally we show $u^r_\i=u^\th_\i=0$. Actually in the cylindrical coordinates, we have the following fact
\be
|\nabla u|^2=|\nabla u^r|^2+|\nabla u^\th|^2+|\nabla u^z|^2+\left|\frac{u^r}{r}\right|^2+\left|\frac{u^\th}{r}\right|^2.\nn
\ee
This means, according to the \eqref{gdc0}, we have
\be\label{2.22}
\int_{-\infty}^\infty\int_0^\infty\left|\frac{u^r}{r}\right|^qrdrdz+\int_{-\infty}^\infty\int_0^\infty\left|\frac{u^\th}{r}\right|^qrdrdz<\infty.
\ee
However, this must be false provided $u^r_\infty$ or $u^\th_\infty$ is non-zero, since we have just proved $u^r(r,z)\to u^r_\infty$ uniformly with respect to $z\in\mathbb{R}$. Therefore if $u^r_\infty\neq 0$, it follows that there exists an $r_0>0$ such that for any $r\geq r_0$,
\[
|u^r(r,z)|\geq\frac{|u^r_\infty|}{2}>0.
\]
This leads to a paradox to \eqref{2.22} since
\be
\int_{-\infty}^\infty\int_0^\infty\left|\frac{u^r}{r}\right|^qrdrdz\geq\int_{-\infty}^\infty\int_{r_0}^\infty\left|\frac{u^r_\infty}{2r}\right|^qrdrdz=\infty.\nn
\ee
The situation of $u^\th$ is similar.
\end{proof}

\subsection{Case $\boldsymbol{q\geq 3}$}
\begin{proposition}
Under the same conditions as those in \textbf{Lemma \ref{DOE}} with $q\geq3$, the following growing estimates of $u$ hold:

if $q=3$, for $r>r_0>0$,
\be\label{uqe4}
|u(r,z)-u(r_0,z)|\leq\,C\log\lt(\f{r}{r_0}+e\rt);
\ee

if $q>3$, for $r>r_0\geq0$,
\be\label{uqb5}
|u(r,z)-u(r_0,z)|\leq\,Cr^{1-3/q},
\ee
where $C$ is independent of $r_0$, $r$ and $z$.

\end{proposition}
\begin{proof}
There exists an $n_0\in\mathbb{N}$ such that $2^{n_0}\leq \f{r}{r_0}\leq2^{n_0+1}$ (note that if $r_0=0$, we let $n_0=+\i$).  Then we iterate the estimate \eqref{OSCEST1} to get the claimed growth of $u$. Here are the details:

If $q=3$ and $r_0>0$,
\be
\begin{split}
|u(r,z)-u(r_0,z)|&\leq|u(r,z)-u(2^{n_0+1}r_0,z)|+\sum_{n=0}^{n_0}\left|u(2^{n+1}r_0,z)-u(2^{n}r_0,z)\right|\\
&\leq\sum_{n=0}^{n_0}\mathop{osc}_{\mathcal{D}_{2^{n}r_0, z}} u\leq C(n_0+1)\leq C\log \lt(\f{r}{r_0}+e\rt),\\
\end{split}\nn
\ee
which indicates \eqref{uqe4}.

Meanwhile, if $q>3$ and $r_0\geq 0$,
\bes
\begin{split}
|u(r,z)-u(r_0,z)|&\leq\sum_{n=0}^{n_0-1}|u(r/2^n,z)-u(r/2^{n+1},z)|+|u(r/2^{n_0},z)-u(r_0,z)|\\
&\leq\sum_{n=0}^\infty\mathop{osc}_{\mathcal{D}_{r/2^n,z}} u\leq C\sum_{n=0}^\infty\left(\frac{r}{2^n}\right)^{1-3/q}\leq Cr^{1-q/3},
\end{split}
\ees
which indicates \eqref{uqb5}.
\end{proof}

\section{Asymptotic behavior of the vorticity: proof of Theorem \ref{thm2}}\label{SEC3}

In 3D Euclidian space (in cylindrical coordinates), for $i\in\bN$, we denote $\mathcal{C}^i_\lambda$ and its related 2D domain $\mathcal{E}^i_\lambda$ by
\bes
\mathcal{C}^i_\lambda:=\left\{(r,\th,z):(1-\f{1}{2^{i+1}})\lambda\leq r\leq (1+\f{1}{2^{i+1}})\lambda,\,0\leq\th< 2\pi,\,|z|\leq \f{1}{2^{i+1}}\lambda\right\},
\ees
and
\bes
\mathcal{E}^i_\lambda:=\left\{(r,z):(1-\f{1}{2^{i+1}})\lambda\leq r\leq (1+\f{1}{2^{i+1}})\lambda,\,|z|\leq \f{1}{2^{i+1}}\lambda\right\},
\ees
respectively. When $\lambda=1$, we write $\mathcal{C}^i$ and $\mathcal{E}^i$ instead of $\mathcal{C}^i_1$ and $\mathcal{E}^i_1$. Before the proof of Theorem \ref{thm2}, the following \emph{Brezis-Wainger} inequality is frequently used.

\begin{lemma}\label{BGI}
Let $\Omega$ be a bounded Lipschitz domain in $\mathbb{R}^2$ and $f\in H^1(\O)\cap W^{1, p}(\O)$ for $p>2$. Then we have
\be\label{BGI1}
\|f\|_{L^\infty(\Omega)}\leq C(\Omega)\log^{1/2}\left(e+\|f\|_{W^{1,p}(\Omega)}\right),
\ee
for every $f\in H^1(\Omega)\cap W^{1, p}(\O)$ with $\|f\|_{H^1(\Omega)}\leq 1$.
\end{lemma}
We refer readers to \cite[Theorem 1]{Brezis-Wainger1980} for details. Although the proof there takes cake of the full space domain $\mathbb{R}^2$, by standard extension arguments, Theorem 1 of \cite{Brezis-Wainger1980} will also hold for the space $H^1(\O)\cap W^{1, p}(\O)$ since a bounded extension operator
\bes
E : H^1(\O)\rightarrow H^1(\bR^2);\q W^{1, p}(\O)\rightarrow W^{1, p}(\bR^2),
\ees
which is a right inverse of the (pointwise) restriction operation, exists if $\p\Omega$ is Lipschitz continuous.

\begin{remark}
For the convenience of the following proof, we will apply
\be\label{BGI2}
\|f\|_{L^\infty(\Omega)}\leq C(\Omega)\left(1+\|f\|_{H^1(\Omega)}\right)\log^{1/2}\left(e+\|\na f\|_{L^3(\Omega)}\right),
\ee
which has no more restriction on the size of $\|f\|_{H^1(\Omega)}$, instead of Lemma \ref{BGI}. Proof of \eqref{BGI2} is to apply \eqref{BGI1} with $f/\|f\|_{H^1(\Omega)}$ and Sobolev embedding $\|f\|_{L^3(\O)}\leq  C(\O)\|f\|_{H^1(\Omega)}$. We omit the details.
\end{remark}

Pick a fixed point $x_0=(\la,0,0)$ for large $\la>0$ in the cylindrical coordinates. Consider the scaled solution
\be
\t{u}(\t{x})=\la u(\la \t{x})=\la u(x),\q \t{w}(\t{x})=\la^2 w(\la \t{x})=\la^2 w(x),\nn
\ee
where $\t{x}=\f{x}{\la}$. Using the scaling-invariant property of the solution of the Navier-Stokes equations, $\t{u}(\t{x}),\ \t{w}(\t{x})$ is also solutions of \eqref{Main} and \eqref{VEQ}. Now we consider $\t{u}(\t{x}),\ \t{w}(\t{x})$ in the domain ${\mathcal{C}^0}$  which correspond to $u(x), w(x)$ in the domain ${\mathcal{C}}^0_\la$. For simplification of notation, we drop the ``$\sim$"  on the scaled solution for a while when computations take place under the scaled sense.

Let $\varphi(x)$ be a cut-off function which satisfies
\[\text{supp}\,\varphi\subset{\mathcal{C}^0}\quad\text{ and }\quad\varphi(x)\equiv 1,\quad\forall\ x \in\mathcal{C}^1,\]
such that derivatives of $\varphi$ up to the second order are bounded. Here goes the proof of Theorem \ref{thm2}.
\subsection{Decay estimate of $\boldsymbol{w^\th}$}\label{SEC3.1}

We test the vorticity equations $\eqref{VEQ}_2$ by $w^\th|w^\th|^{q-2}\varphi^q$, then it follows that

\be\label{wth0}
\begin{split}
&\q\int_{{\mathcal{C}^0}}w^\th|w^\th|^{q-2}\varphi^q\left(\Delta-\frac{1}{r^2}\right)w^\th dx\\
&=\int_{{\mathcal{C}^0}}\left[b\cdot\nabla w^\th\cdot w^\th|w^\th|^{q-2}\varphi^q-\frac{u^r}{r}\left|w^\th\right|^q\varphi^q+2\frac{w^r}{r}u^\th w^\th|w^\th|^{q-2}\varphi^q\right]dx.
\end{split}
\ee
Using integration by parts and H\"{o}lder inequality, the above equality leads to the following inequality

\be\label{wth}
\begin{split}
&\q\int_{{\mathcal{C}^0}}\left|\nabla(w^\th\varphi)\right|^2\cdot|w^\th\varphi|^{q-2}dx
+\int_{{\mathcal{C}^0}}\frac{1}{r^2}\left|(w^\th\varphi)\right|^qdx\\
&\ls\left(1+\left\|(u^r,u^\th,u^z)\right\|_{L^\infty({\mathcal{C}^0})}\right)\left(\int_{{\mathcal{C}^0}}\left|w^\th\right|^q dx
+\int_{{\mathcal{C}^0}}\left|w^r\right|^qdx\right).
\end{split}
\ee
When we estimate \eqref{wth0}, we use the fact $r\thickapprox 1$ in ${\mathcal{C}^0}$. By the definition of the cut-off function $\varphi$, one finds \eqref{wth} lead to,
\be\label{wth1}
\int_{\mathcal{C}^1}\left|\nabla\left(w^\th\right)^{q/2}\right|^2dx\lesssim \left(1+\left\|(u^r,u^\th,u^z)\right\|_{L^\infty({\mathcal{C}^0})}\right)\left(\int_{\mathcal{C}^0}\left|w^\th\right|^qdx
+\int_{{\mathcal{C}^0}}\left|w^r\right|^qdx\right).
\ee
Since $w^\th$ depends only on $r,z$ and in ${\mathcal{C}^0}$, $r\thickapprox 1$, then we have
\be\label{2.62e}
\begin{aligned}
&\q\int_{\mathcal{E}^1}\left|\bar{\nabla}\left(w^\th\right)^{q/2}\right|^2drdz\\
&\lesssim \left(1+\left\|(u^r,u^\th,u^z)\right\|_{L^\infty({\mathcal{E}^0})}\right)\left(\int_{{\mathcal{E}^0}}\left|w^\th\right|^qdrdz
+\int_{{\mathcal{E}^0}}\left|w^r\right|^qdrdz\right),
\end{aligned}
\ee
where $\bar{\na}=(\p_r, \p_z)$. By the \emph{Brezis-Wainger}-type inequality \eqref{BGI2}, one derives,
\be\label{wthbg}
\left\|\left(w^\th\right)^{q/2}\right\|_{L^\infty(\mathcal{E}^1)}\lesssim\left(1+\left\|\left(w^\th\right)^{q/2}\right\|_{H^1(\mathcal{E}^1)}\right)
\log^{1/2}\left(e+\left\|\na\left(w^\th\right)^{q/2}\right\|_{L^3(\mathcal{E}^1)}\right)
\ee
by choosing $f=\left(w^\th\right)^{q/2}$. Now using \eqref{2.62e} and \eqref{wthbg} and going back to the 3-dimensional domain $\mathcal{C}^1$, we have
\be
\begin{split}
\left\|\left(w^\th\right)^{q/2}\right\|_{L^\infty(\mathcal{C}^1)}\lesssim&\left(1+\left\|(u^r,u^\th,u^z)\right\|^{1/2}_{L^\infty({\mathcal{C}^0})}\right)
\left(\int_{{\mathcal{C}^0}}\left|\na u\right|^qdx\right)^{1/2}\\
&\cdot\log^{1/2}\left(e+\left\|\na\left(w^\th\right)^{q/2}\right\|_{L^3(\mathcal{C}^1)}\right).
\end{split}\nn
\ee
Now we need to bound $\left\|\na \left(w^\th\right)^{q/2}\right\|_{L^3(\mathcal{C}^1)}$. Actually we will see in \eqref{erea1},
\be\label{ere1}
\left\|\nabla\left(w^\th\right)^{q/2}\right\|_{L^3(\mathcal{C}^1)}\leq \mathcal{A}\lt(\|u\|_{L^\i({\mathcal{C}^0})},\|\nabla u\|_{L^q({\mathcal{C}^0})} \rt),
\ee
where $\mathcal{A}$ is a positive power function depending on $\|u\|_{L^\i({\mathcal{C}^0})}$ and $\|\nabla u\|_{L^q({\mathcal{C}^0})}$ whose explicit representation is not important for us. After scaling back to the domains with  ``$\lambda-$size" for $\lambda>>1$, it can only grow at most as a polynormial order of $\la$ at the far field. Since it appears in a ``$\log$" function, we need not to calculate the exact order. The calculation of \eqref{ere1} is presented in Appendix \ref{appa}.

Now we take back the ``$\sim$" to the scaled solution and apply \eqref{ere1}, then we have
\be
\begin{split}
\left\|\left(\t{w}^\th\right)^{q/2}\right\|_{L^\infty(\mathcal{C}^1)}
\lesssim&\left(1+\left\|(\t{u}^r,\t{u}^\th,\t{u}^z)\right\|^{1/2}_{L^\infty({\mathcal{C}^0})}\right)\left(\int_{{\mathcal{C}^0}}
\left|\t{\na}\t{u}\right|^qd\t{x}\right)^{1/2}\\
&\cdot\log^{1/2}\left(e+\mathcal{A}\lt(\|\t{u}\|_{L^\i({\mathcal{C}^0})},\|\t{\nabla}\t{u}\|_{L^q({\mathcal{C}^0})} \rt)\right).
\end{split}\nn
\ee
If we scale back to the domains with  ``$\lambda-$size" for $\lambda>>1$, then we have
\begin{footnotesize}
\be\label{3.14}
\begin{split}
\lambda^q\left\|\left({w}^\th\right)^{q/2}\right\|_{L^\infty(\mathcal{C}^1_\lambda)}
\lesssim&\left(1+\lambda^{1/2}\left\|({u}^r,{u}^\th,{u}^z)\right\|^{1/2}_{L^\infty(\mathcal{C}^0_\lambda)}\right)
\cdot\lambda^{q-3/2}\left(\int_{\mathcal{C}^0_\lambda}\left|\na u\right|^qdx\right)^{1/2}\\
&\cdot\log^{1/2}\left(e+\la^M\mathcal{A}\lt(\|{u}\|_{L^\i({\mathcal{C}^0_\la})},\|{\nabla}{u}\|_{L^q({\mathcal{C}^0_\la})} \rt)\right),
\end{split}\nn
\ee\end{footnotesize}
where $M$ is the scaling power of $\mathcal{A}$, whose exact value is not important for us since it appears inside a ``$log$". Therefore $w^\th$ decays as
\be\label{wth2}
\left\|w^\th\right\|_{L^\infty(\mathcal{C}^1_\lambda)}\lesssim\left(1+\lambda^{1/2}
\left\|({u}^r,{u}^\th,{u}^z)\right\|^{1/2}_{L^\infty(\mathcal{C}^0_\lambda)}\right)^{2/q}\lambda^{-3/q}\left(\log\lambda\right)^{1/q}.
\ee
\noindent \textbf{Case I:} Under Assumption \ref{assum} and in Case I of Theorem \ref{thm2}, by using \eqref{uqs3}, \eqref{uqe3} and \eqref{uqb3}, we see that for $\la$ large,
\be\label{uqb32}
\left\|({u}^r,{u}^\th,{u}^z)\right\|_{L^\infty(\mathcal{C}^0_\lambda)}\leq\lt\{
\begin{aligned}
&C\la^{1-3/q},\,\,\,\,\q\quad\quad\text{for }q\in(2,3);\\
&C({r_0})\log\la,\,\,\,\q\quad\text{for }q=3;\\
&C({r_0})\la^{1-3/q},\q\quad\text{for }q\in(3,\infty).
\end{aligned}
\rt.
\ee
Inserting \eqref{uqb32} into \eqref{wth2}, we can get
\be
\left\|w^\th\right\|_{L^\infty(\mathcal{C}^1_\lambda)}\lesssim\lt\{
\begin{aligned}
&\la^{-1/q-3/q^2}\left(\log\la\right)^{1/q},\,\,\,\,\,\,\,\,\quad\text{for }q\in(2,3)\cup(3,\infty);\\
&\la^{-2/3}(\log\la)^{2/3},\,\,\,\q\q\q\quad\text{for }q=3,\\
\end{aligned}
\rt.\nn
\ee
which indicates the estimate of $w^\th$ in \eqref{westimate1}.

\noindent \textbf{Case II:} In Case II of Theorem \ref{thm2}, by using \eqref{uqs3}, we see that for $\la$ large

\be\label{uqs33}
\left\|({u}^r,{u}^\th,{u}^z)\right\|_{L^\infty(\mathcal{C}^0_\lambda)}\leq C.
\ee
Inserting \eqref{uqs33} into \eqref{wth2}, we can get
\be
\left\|w^\th\right\|_{L^\infty(\mathcal{C}^1_{\lambda})}\lesssim\lambda^{-2/q}\left(\log\lambda\right)^{1/q},\nn
\ee
which indicates the estimate of $w^\th$ in \eqref{westimate2}.

\subsection{Decay estimates of $\boldsymbol{w^r}$ and $\boldsymbol{w^z}$}
Decay estimates of $w^r$ and $w^z$ are much more involved and can not be as good as that for $w^\th$. We will use an iteration technique to get a $(\f{1}{q-1})^-$ order decay for $q\geq 3$. Then for $q\in (2,3)$, actually this order can be improved by using the decay of $w^\th$, $\na w^\th$ and a pointwise estimate Lemma concerning Calderon-Zygmund operator given in \cite{CPZ2018}.

First we perform some energy estimates for the scaled vorticity $\t{w}^r$ and $\t{w}^z$ (still denote them by ${w}^r$ and ${w}^z$). We test the vorticity equations $\eqref{VEQ}_1$ and $\eqref{VEQ}_3$ by $w^r|w^r|^{q-2}\varphi^q$ and $w^z|w^z|^{q-2}\varphi^q$ respectively, then it follows that
\be
\begin{split}
&\q\int_{{\mathcal{C}^0}}w^r|w^r|^{q-2}\varphi^q\left(\Delta-\frac{1}{r^2}\right)w^rdx\\
&=\int_{{\mathcal{C}^0}}\left[b\cdot\nabla w^r\cdot w^r|w^r|^{q-2}\varphi^q-(w^r\p_r+w^z\p_z)u^r\cdot w^r|w^r|^{q-2}\varphi^q\right]dx,
\end{split}\nn
\ee
\be
\begin{split}
&\q\int_{{\mathcal{C}^0}}w^z|w^z|^{q-2}\varphi^q\Delta w^zdx\\
&=\int_{{\mathcal{C}^0}}\left[b\cdot\nabla w^z\cdot w^z|w^z|^{q-2}\varphi^q-(w^r\p_r+w^z\p_z)u^z\cdot w^z|w^z|^{q-2}\varphi^q\right]dx.
\end{split}\nn
\ee
Using integration by parts, Young inequality and H\"{o}lder inequality, the above equations of $w^r$ and $w^z$ lead to the following 3 groups of inequalities. First,
\bes
\bali
&\q\int_{{\mathcal{C}^0}}\left|\nabla(w^r\varphi)^{q/2}\right|^2dx+\int_{{\mathcal{C}^0}}\frac{1}{r^2}\left|(w^r\varphi)\right|^qdx\\
&\lesssim \left(1+\left\|(u^r,u^z)\right\|_{L^\infty({\mathcal{C}^0})}+\left\|(\nabla u^r,\nabla u^z)\right\|_{L^\infty({\mathcal{C}^0})}\right)\left(\int_{{\mathcal{C}^0}}\left|w^r\right|^qdx+\int_{{\mathcal{C}^0}}\left|w^z\right|^qdx\right),
\eali
\ees
\bes
\bali
&\q\int_{{\mathcal{C}^0}}\left|\nabla(w^z\varphi)^{q/2}\right|^2dx\\
&\lesssim \left(1+\left\|(u^r,u^z)\right\|_{L^\infty({\mathcal{C}^0})}+\left\|(\nabla u^r,\nabla u^z)\right\|_{L^\infty({\mathcal{C}^0})}\right)\left(\int_{{\mathcal{C}^0}}\left|w^r\right|^qdx+\int_{{\mathcal{C}^0}}\left|w^z\right|^qdx\right).
\eali
\ees
Second,
\bes
\bali
&\q\int_{{\mathcal{C}^0}}\left|\nabla(w^r\varphi)^{q/2}\right|^2dx+\int_{{\mathcal{C}^0}}\frac{1}{r^2}\left|(w^r\varphi)\right|^qdx\\
&\lesssim \left(1+\left\|(u^r,u^z)\right\|_{L^\infty({\mathcal{C}^0})}+\left\|(w^r,w^z)\right\|_{L^\infty({\mathcal{C}^0})}\right)\int_{{\mathcal{C}^0}}\left|\na u\right|^qdx,
\eali
\ees
\bes
\bali
&\q\int_{{\mathcal{C}^0}}\left|\nabla(w^z\varphi)^{q/2}\right|^2dx\\
&\lesssim \left(1+\left\|(u^r,u^z)\right\|_{L^\infty({\mathcal{C}^0})}+\left\|(w^r,w^z)\right\|_{L^\infty({\mathcal{C}^0})}\right)\int_{{\mathcal{C}^0}}
\left|\nabla u \right|^qdx.
\eali
\ees
Third,
\bes
\int_{{\mathcal{C}^0}}\left|\nabla(w^r\varphi)^{q/2}\right|^2dx+\int_{{\mathcal{C}^0}}\frac{1}{r^2}\left|(w^r\varphi)\right|^qdx\lesssim \left(1+\left\|(u^r,u^z)\right\|_{L^\infty({\mathcal{C}^0})}\right)^2 \int_{{\mathcal{C}^0}}\left|\na u\right|^qdx,
\ees
\bes
\int_{{\mathcal{C}^0}}\left|\nabla(w^z\varphi)^{q/2}\right|^2dx\lesssim \left(1+\left\|(u^r,u^z)\right\|_{L^\infty({\mathcal{C}^0})}\right)^2 \int_{{\mathcal{C}^0}}\left|\na u\right|^qdx.
\ees
Adding each group above together and noting $\varphi=1$ in $\mathcal{C}^1$, we get
{\begin{footnotesize}
\be\label{3.1777}
\int_{\mathcal{C}^1}\left|\nabla\left(w^r, w^z\right)^{q/2}\right|^2dx\ls\left\{
\begin{aligned}
&\left(1+\left\|(u^r,u^z)\right\|_{L^\infty({\mathcal{C}^0})}+\left\|(\nabla u^r,\nabla u^z)\right\|_{L^\infty({\mathcal{C}^0})}\right)\int_{\mathcal{C}^0}|\nabla u|^qdx;\\[2mm]
&\left(1+\left\|(u^r,u^z)\right\|_{L^\infty({\mathcal{C}^0})}+\left\|(w^r,w^z)\right\|_{L^\infty({\mathcal{C}^0})}\right)
\int_{{\mathcal{C}^0}}\left|\na u\right|^qdx;\\[2mm]
&\left(1+\left\|(u^r,u^z)\right\|^2_{L^\infty({\mathcal{C}^0})}\right)\int_{{\mathcal{C}^0}}\left|\na u\right|^qdx.
\end{aligned}
\rt.
\ee
\end{footnotesize}}
Since $w^r, w^z$ depends only on $r,z$ and in ${\mathcal{C}^0}$, $r\thickapprox 1$, then we have related 2D estimate of \eqref{3.1777}
{\begin{footnotesize}
\be\label{3.1888}
\int_{\mathcal{E}^1}\left|\bar{\nabla}\left(w^r, w^z\right)^{q/2}\right|^2drdz\ls\left\{
\begin{aligned}
&\left(1+\left\|(u^r,u^z)\right\|_{L^\infty({\mathcal{E}^0})}+\left\|(\bar\nabla u^r,\bar\nabla u^z)\right\|_{L^\infty({\mathcal{E}^0})}\right)\int_{{\mathcal{E}^0}}\left|\bar{\na} u\right|^qdrdz;\\[2mm]
&\left(1+\left\|u^r,u^z\right\|_{L^\infty({\mathcal{E}^0})}+\left\|(w^r,w^z)\right\|_{L^\infty({\mathcal{E}^0})}\right)\int_{{\mathcal{E}^0}}\left|\bar{\na} u\right|^qdrdz;\\[2mm]
&\left(1+\left\|(u^r,u^z)\right\|^2_{L^\infty({\mathcal{E}^0})}\right)\int_{{\mathcal{E}^0}}\left|\bar{\na} u\right|^qdrdz.
\end{aligned}
\rt.
\ee
\end{footnotesize}
}
where $\bar{\na}=(\p_r,\p_z)$. By \eqref{BGI2}, one derives,
{\begin{footnotesize}
\be\label{wrwz1}
\bali
\left\|\left(w^r,w^z\right)^{q/2}\right\|_{L^\infty(\mathcal{E}^1)}\lesssim\left(1+\left\|\left(w^r,w^z\right)^{q/2}\right\|_{H^1(\mathcal{E}^1)}\right)
\log^{1/2}\left(e+\left\|\nabla\left(w^r,w^z\right)^{q/2}\right\|_{L^3(\mathcal{E}^1)}\right).
\eali
\ee
\end{footnotesize}}
Now using $\eqref{3.1888}$, \eqref{wrwz1} and going back to the 3-dimensional domains $\mathcal{C}^0,\,\mathcal{C}^1$, we have
{\begin{footnotesize}
\be\label{3.2000}
\left\|\left(w^r,w^z\right)^{q/2}\right\|_{L^\infty(\mathcal{C}^1)}\ls\left\{
\begin{aligned}
&\left(1+\left\|(u^r,u^z)\right\|^{1/2}_{L^\infty({\mathcal{C}^0})}+\left\|(\nabla u^r,\nabla u^z)\right\|^{1/2}_{L^\infty({\mathcal{C}^0})}\right)\times\\
&\q\left(\int_{{\mathcal{C}^0}}\left|\na u\right|^qdx\right)^{1/2}\log^{1/2}\left(e+\left\|\na \left(w^r,w^z\right)^{q/2}\right\|_{L^3({\mathcal{C}^1})}\right);\\[4mm]
&\left(1+\left\|(u^r,u^z)\right\|^{1/2}_{L^\infty({\mathcal{C}^0})}+\left\|(w^r,w^z)\right\|^{1/2}_{L^\infty({\mathcal{C}^0})}\right)\times
\\
&\q\left(\int_{{\mathcal{C}^0}}\left|\na u\right|^qdx\right)^{1/2}\log^{1/2}\left(e+\left\|\na \left(w^r,w^z\right)^{q/2}\right\|_{L^3({\mathcal{C}^1})}\right);\\[4mm]
&\left(1+\left\|(u^r,u^z)\right\|_{L^\infty({\mathcal{C}^0})}\right)\times\\
&\q\left(\int_{{\mathcal{C}^0}}\left|\na u\right|^qdx\right)^{1/2}\log^{1/2}\left(e+\left\|\na \left(w^r,w^z\right)^{q/2}\right\|_{L^3({\mathcal{C}^1})}\right).
\end{aligned}
\rt.\nn
\ee
\end{footnotesize}
}
Now we take back the ``$\sim$" to the scaled solution, it follows that
{\begin{footnotesize}
\be\label{3.2111}
\left\|\left(\t{w}^r,\t{w}^z\right)^{q/2}\right\|_{L^\infty(\mathcal{C}^1)}\ls\left\{
\begin{array}{l}
\left(1+\left\|(\t{u}^r,\t{u}^z)\right\|^{1/2}_{L^\infty({\mathcal{C}^0})}+\left\|(\t{\nabla} \t{u}^r,\t{\nabla} \t{u}^z)\right\|^{1/2}_{L^\infty({\mathcal{C}^0})}\right)\times\\
\q\left(\int_{{\mathcal{C}^0}}\left|\t{\na} \t{u}\right|^qd\t{x}\right)^{1/2}\log^{1/2}\left(e+\left\|\t{\na} \left(\t{w}^r,\t{w}^z\right)^{q/2}\right\|_{L^3({\mathcal{C}^1})}\right);\\[4mm]
\left(1+\left\|(\t{u}^r,\t{u}^z)\right\|^{1/2}_{L^\infty({\mathcal{C}^0})}+\left\|( \t{w}^r, \t{w}^z)\right\|^{1/2}_{L^\infty({\mathcal{C}^0})}\right)\times\\
\q\left(\int_{{\mathcal{C}^0}}\left|\t{\na} \t{u}\right|^qd\t{x}\right)^{1/2}\log^{1/2}\left(e+\left\|\t{\na} \left(\t{w}^r,\t{w}^z\right)^{q/2}\right\|_{L^3({\mathcal{C}^1})}\right);\\[4mm]
\left(1+\left\|(\t{u}^r,\t{u}^z)\right\|_{L^\infty({\mathcal{C}^0})}\right)\times\\
\q\left(\int_{{\mathcal{C}^0}}\left|\t{\na} \t{u}\right|^q d\t{x}\right)^{1/2}\log^{1/2}\left(e+\left\|\t{\na} \left(\t{w}^r,\t{w}^z\right)^{q/2}\right\|_{L^3({\mathcal{C}^1})}\right).
\end{array}
\rt.\nn
\ee
\end{footnotesize}
}
If we scale back to the domains with ``$\lambda-$size" for $\lambda>>1$, then we have achieved
\begin{footnotesize}
\be\label{3.2222}
\la^q\left\|\left({w}^r,{w}^z\right)^{q/2}\right\|_{L^\infty(\mathcal{C}^1_\la)}\ls\left\{
\begin{array}{l}
\left(1+\la^{1/2}\left\|({u}^r,{u}^z)\right\|^{1/2}_{L^\infty({\mathcal{C}}^0_\la)}+\la\left\|(\nabla {u}^r,\nabla {u}^z)\right\|^{1/2}_{L^\infty(\mathcal{C}^0_\la)}\right)\la^{q-3/2}\times\\
\q\left(\int_{\mathcal{C}^0_\la}\left|\na u\right|^qdxdx\right)^{1/2}\log^{1/2}\left(e+\la^{q}\left\|\na \left({w}^r,{w}^z\right)^{q/2}\right\|_{L^3({\mathcal{C}}^1_\la)}\right);\\[4mm]
\left(1+\la^{1/2}\left\|({u}^r,{u}^z)\right\|^{1/2}_{L^\infty({\mathcal{C}}^0_\la)}+\la\left\|(w^r,w^z)\right\|^{1/2}_{L^\infty(\mathcal{C}^0_\la)}\right)\la^{q-3/2}\times\\
\q\left(\int_{\mathcal{C}^0_\la}\left|\na u\right|^qdxdx\right)^{1/2}\log^{1/2}\left(e+\la^{q}\left\|\na \left({w}^r,{w}^z\right)^{q/2}\right\|_{L^3({\mathcal{C}}^1_\la)}\right);\\[4mm]
\left(1+\la\left\|({u}^r,{u}^z)\right\|_{L^\infty({\mathcal{C}}^0_\la)}\right)\la^{q-3/2}\times\\
\q \left(\int_{\mathcal{C}^0_\la}\left|\na u\right|^qdxdx\right)^{1/2}\log^{1/2}\left(e+\la^{q}\left\|\na \left({w}^r,{w}^z\right)^{q/2}\right\|_{L^3({\mathcal{C}}^1_\la)}\right).
\end{array}
\rt.\nn
\ee
\end{footnotesize}
Similarly as in Section \ref{SEC3.1}, $(w^r, w^z)$ decays as
{\begin{footnotesize}
\be\label{3.2333}
\left\|(w^r,w^z)\right\|_{L^\infty(\mathcal{C}^1_\lambda)}\ls\left\{
\begin{array}{l}
\left(\la^{1/2}\left\|({u}^r,{u}^z)\right\|^{1/2}_{L^\infty(\mathcal{C}^0_\la)}+\la\left\|(\nabla {u}^r,\nabla {u}^z)\right\|^{1/2}_{L^\infty(\mathcal{C}^0_\la)}\right)^{2/q}\lambda^{-3/q}\left(\log\lambda\right)^{1/q};\\[2mm]
\left(\la^{1/2}\left\|({u}^r,{u}^z)\right\|^{1/2}_{L^\infty(\mathcal{C}^0_\la)}
+\la\left\|(w^r,w^z)\right\|^{1/2}_{L^\infty(\mathcal{C}^0_\la)}\right)^{2/q}\lambda^{-3/q}\left(\log\lambda\right)^{1/q};\\[2mm]
\left(\la\left\|({u}^r,{u}^z)\right\|_{L^\infty(\mathcal{C}^0_\la)}\right)^{2/q}\lambda^{-3/q}\left(\log\lambda\right)^{1/q}.
\end{array}
\rt.
\ee
\end{footnotesize}
}
\begin{remark}
 The reason why decay estimates of $w^r$ and $w^z$ are weaker than that of $w^\th$ is due to $\eqref{3.2333}_3$ where $\la\left\|({u}^r,{u}^z)\right\|_{L^\infty(\mathcal{C}^0_\la)}$ here is replaced by $\la^{1/2}\left\|({u}^r,{u}^z)\right\|^{1/2}_{L^\infty(\mathcal{C}^0_\la)}$ in the decay estimate inequality \eqref{wth2} of $w^\th$. Now we use $\eqref{3.2333}_3$ to provide an iteration initial decay and then use $\eqref{3.2333}_2$ to iterate the decay estimates of $w^r$ and $w^z$ on $\mathcal{C}^i_{\la}$ for $1\leq i\in\bN$. After finite times' iterations, we can achieve that for $q\geq 3$,
 \bes
 |(w^r,w^z)|=O(r^{-(\f{1}{q-1})^{-}}), \q \text{as }\ r\rightarrow +\i.
 \ees
 While when $q\in (2,3)$, the above decay estimates can be improved. By using the decay estimates of $w^\th$, $\na w^\th$ and Lemma \ref{Lem3.3} below, we get a decay estimate for $\na u^r$ and $\na u^z$, then  inserting this estimate into $\eqref{3.2333}_1$, we can get an improved estimate of $w^r$ and $w^z$.
\end{remark}




\subsubsection{Decay by iteration for $\boldsymbol{q\geq 3}$ }
From \eqref{uqb32} and $\eqref{3.2333}_3$, we have in {\bf Case I} of Theorem \ref{thm2}, for $q\geq 3$
\be\label{3.1555555}
\left\|(w^r,w^z)\right\|_{L^\infty(\mathcal{C}^1_\lambda)}\leq C\la^{\lt(\f{1}{q}-\f{6}{q^2}\rt)}\lt(\log\la\rt)^{\f{3}{q}}.
\ee
When $q>6$, we see that the order $\f{1}{q}-\f{6}{q^2}$ is positive and this estimate is very bad.


Next we use iteration to improve the decay order. From the above estimate, we first have
\be\label{wrwzinitial}
\left\|(w^r,w^z)\right\|_{L^\infty(\mathcal{C}^1_\lambda)}\leq C\la^{(1/q-6/q^2)^+} .
\ee
 Here and below we denote by $\al^+$ a constant which is larger but close to $\al$. Using $\eqref{3.2333}_2$, we have

\bes
\bali
&\q \left\|(w^r,w^z)\right\|_{L^\infty(\mathcal{C}^1_\lambda)}\\
&\lesssim\lt[\left(\la^{1/2}\left\|({u}^r,{u}^z)\right\|^{1/2}_{L^\infty(\mathcal{C}^0_\la)}\rt)^{2/q}+\lt(
\la\left\|(w^r,w^z)\right\|^{1/2}_{L^\infty(\mathcal{C}^0_\la)}\right)^{2/q}\rt]\lambda^{-3/q}\left(\log\lambda\right)^{1/q}\\
&\leq C_1\lt\{ \la^{-(1/q+3/q^2)}\left(\log\lambda\right)^{2/q}+\la^{-1/q}\left(\log\lambda\right)^{1/q}\left\|(w^r,w^z)\right\|_{L^\infty(\mathcal{C}^0_\la)}^{1/q}\rt\}.
\eali
\ees
Actually for any $i\in \bN$, we have a constant $C_i$, which is independent of $\lambda$ and will go to infinity as $i\rightarrow+\i$, such that
\be\label{3.16555}
\bali
&\q \left\|(w^r,w^z)\right\|_{L^\infty(\mathcal{C}^i_\lambda)}\\
&\leq C_i\lt\{ \la^{-(1/q+3/q^2)}\left(\log\lambda\right)^{2/q}+\la^{-1/q}\left(\log\lambda\right)^{1/q}\left\|(w^r,w^z)\right\|_{L^\infty(\mathcal{C}^{i-1}_\la)}^{1/q}\rt\}.
\eali
\ee
 Since $\f{1}{q}+\f{3}{q^2}>\f{1}{q-1}>\f{1}{q}-\f{6}{q^2}$ for all $q\geq 3$, if we start with the estimate \eqref{3.1555555} and iterate over \eqref{3.16555}, the decay of $(w^r,w^z)$ will always be refined each time. Thus for some $C_i>0$, we have
\be
\bali
\left\|(w^r,w^z)\right\|_{L^\infty(\mathcal{C}^i_\lambda)}\leq C_i\lt\{\left(\f{\log\lambda}{\la}\right)^{1/q}\left\|(w^r,w^z)\right\|_{L^\infty(\mathcal{C}^{i-1}_\la)}^{1/q}\rt\},\quad\forall i\in\mathbb{N}.
\eali\nn
\ee
After $n$ times iteration, we have
\bes
\bali
\left\|(w^r,w^z)\right\|_{L^\infty(\mathcal{C}^n_\lambda)}&\leq \prod_{i=1}^nC_i\lt\{\left(\f{\log\lambda}{\la}\right)^{\sum\limits^{n-1}_{k=1}(1/q)^k}\left\|(w^r,w^z)\right\|_{L^\infty(\mathcal{C}^{1}_\la)}^{(1/q)^{n-1}}\rt\}\\
&\leq \prod_{i=1}^nC_i \left(\f{\log\lambda}{\la}\right)^{\sum\limits^{n-1}_{k=1}(1/q)^k} \lt(\la^{(\f{1}{q}-\f{6}{q^2})^+}\rt)^{(1/q)^{n-1} },
\eali
\ees
where we have used the initial estimate \eqref{wrwzinitial}. By noting that

 \[
\lim_{n\rightarrow+\i} (1/q)^{n-1}(1/q-6/q^2)^+=0
\]
and $\sum\limits^{\i}_{k=1}(1/q)^k=1/(q-1)$, then for sufficiently large $n$, the above iteration indicates that
\be\label{wrwzrev6}
\left\|(w^r,w^z)\right\|_{L^\infty(\mathcal{C}^n_\la)}\les\la^{-\lt(\f{1}{q-1}\rt)^{-}}.
\ee
 Thus we conclude that
\be\label{eq>3}
 |(w^r,w^z)|=O(r^{-(\f{1}{q-1})^{-}}), \q \text{as }\ r\rightarrow +\i.
\ee

\subsubsection{Improved decay for $\boldsymbol{2<q< 3}$}
Now we use $\eqref{3.2333}_1$ to improve the decay estimate for $q\in (2,3)$. To derive the $L^\infty$ estimates of $\nabla u^r$ and $\nabla u^z$, we note that by using Biot-Savart law,
\be
\nabla(u^re_r+u^ze_z)=\nabla(-\Delta)^{-1}\text{curl}(w^\th e_\th),\nn
\ee
which implies
\be
\nabla u^r=\int_{\mathbb{R}^3}K_1(x,y)w^\th(y)dy,\quad\nabla u^z=\int_{\mathbb{R}^3}K_2(x,y)w^\th(y)dy,\nn
\ee
where $K_1$ and $K_2$ are Calderon-Zygmund kernels. The following lemma describes the property of the Calderon-Zygmund kernels act on axially symmetric functions which may help us derive the decay estimates of $\nabla u^r$ and $\nabla u^z$.
\begin{lemma}[c.f. Lemma 3.2 in \cite{CPZ2018}]\label{Lem3.3}
Assume that $K(x,y)$ be a Calderon-Zygmund kernel and $f$ is a smooth axially symmetric function satisfying, for $x=(x',z)\in\mathbb{R}^3$
\[
|f(x)|+|\nabla f(x)|\lesssim\frac{\log^b(e+|x'|)}{(1+|x'|)^a},\quad\text{for }0<a<2,\,\,b>0.
\]
Define $Tf(x):=\int K(x,y)f(y)dy$. Then there exists a constant $C_0$ such that
\be
|Tf(x)|\leq C_0\frac{\log^{b+1}(e+|x'|)}{(1+|x'|)^a}.\nn
\ee
\end{lemma}

\noindent \textbf{Case I:} Under Assumption \ref{assum} and $u^z_\i=0$ in Theorem \ref{thm2}, carrying out the similar estimate as that of $w^\th$ in the $L^2$ framework, we can show that from \eqref{reva1}
\bes
\left\|\na w^\th\right\|_{L^\infty(\mathcal{C}^2_\lambda)}\lesssim
\begin{aligned}
&\la^{-(15/2q-2)}\left(\log\la\right)^{1/2},\,\,\,\,\,\,\,\,\quad\text{for }q\in(2,3).
\end{aligned}
\ees
The details are given in Appendix \ref{appb}.
Applying Lemma \ref{Lem3.3} with $f=w^\th$ and $K(x,y)=K_1(x,y)$ and $K_2(x,y)$, we can get, for large $r$, that
\be\label{duruz}
|\na u^r(r,z)|+|\na u^z(r,z)|\leq\lt\{
\begin{aligned}
&r^{-1/q-3/q^2}\left(\log r \right)^{3/2},\,\,\quad\text{for }q\in(2,\f{13+\s{73}}{8});\\
&r^{-(15/2q-2)}\left(\log r \right)^{3/2},\,\,\quad\text{for }q\in[\f{13+\s{73}}{8},3).
\end{aligned}
\rt.
\ee
Inserting \eqref{uqb32} and \eqref{duruz} into $\eqref{3.2333}_1$, we can get
\bes
\left\|(w^r,w^z)\right\|_{L^\infty(\mathcal{C}^4_\lambda)}\lesssim\lt\{
\begin{aligned}
&\la^{-(1/q+1/q^2+3/q^3)^-},\,\,\quad\text{for }q\in(2,\f{13+\s{73}}{8});\\
&\la^{-(-1/q+15/2q^2)^-},\,\,\quad\text{for }q\in[\f{13+\s{73}}{8},3),
\end{aligned}
\rt.
\ees
which together \eqref{eq>3} indicate the estimates of $w^r$ and $w^z$ in \eqref{westimate1}.

\noindent \textbf{Case II:} In this case, $u$ is bounded, using \eqref{reva1}, it is routine to show that
\be
\left\|\na w^\th\right\|_{L^\infty(\mathcal{C}^2_\lambda)}\lesssim\lambda^{1/2-3/q}\left(\log\lambda\right)^{1/2}.\nn
\ee
Applying Lemma \ref{Lem3.3} with $f=w^\th$ and $K(x,y)=K_1(x,y)$ and $K_2(x,y)$, we can get
\be\label{duruz1}
|\na u^r(r,z)|+|\na u^z(r,z)|\leq  Cr^{-(3/q-1/2)}\left(\log r\right)^{3/2},\q \text{for\ large}\ r.
\ee
Inserting \eqref{uqs33} and \eqref{duruz1} into $\eqref{3.2333}_1$, we can get
\be
\left\|(w^r,w^z)\right\|_{L^\infty(\mathcal{C}^4_\lambda)}\lesssim\lambda^{-(1/2q+3/q^2)^-},\nn
\ee
which indicates the estimates of $w^r$ and $w^z$ in \eqref{westimate2}.

\section*{Acknowledgments} The authors wish to thank Prof. Qi S. Zhang in UC Riverside for his constant encouragement on this topic and the referees for careful reading and very useful comments. Z. Li is supported by Natural Science Foundation of Jiangsu Province (No. BK20200803) and National Natural Science Foundation of China (No. 12001285). X. Pan is supported by National Natural Science Foundation of China (No. 12031006, No. 11801268).

\appendix

\section{Estimates of $\boldsymbol{\|\na w^{q/2}\|_{L^3(\mathcal{C}^1)}}$}\label{appa}
First we derive $L^2$-norm estimates for derivatives of $u$.

\begin{lemma}\label{A1}
 For $i\in\bN$ and $0\leq i\leq 8$, define $\mathfrak{C}_i$ as
\bes
\mathfrak{C}_i:=\left\{(r,\th,z):\f{1}{2}+\f{i}{32}\leq r\leq \f{3}{2}-\f{i}{32},\,0\leq\th< 2\pi,\,|z|\leq \f{1}{2}-\f{i}{32}\right\}.
\ees
We see that $\mathfrak{C}_8=\mathcal{C}^1$ and $\mathfrak{C}_0={\mathcal{C}^0}$. Then
\be\label{eredl2}
 \left\|\na^k u\right\|_{L^2(\mathfrak{C}_6)}  \leq \mathcal{A}\lt(\|u\|_{L^\i({\mathcal{C}^0})},\|\nabla u\|_{L^q({\mathcal{C}^0})} \rt)
\ee
for  $q\geq2$, $k=2,3,4$.
\end{lemma}
\begin{proof}
In the full 3-dimensional space, $w$ satisfies
\be\label{erev}
-\Dl w+u\cdot\na w-w\cdot\na u=0.
\ee
 Let $\varphi_i(x)$ $(1\leq i\leq 8)$ be a cut-off function which satisfies
\[
\text{supp}\,\varphi_i\subset\mathfrak{C}_{i-1}\quad\text{ and }\quad\varphi(x)\equiv 1,\quad\forall x\in\mathfrak{C}_i.
\]
Testing \eqref{erev} by $w\varphi^2_1$ and then using integration by parts indicate that
\be
\bali
\int |\na w|^2\varphi^2_1dx+\int w\na w \na\varphi^2_1dx=\int \lt(w\cdot\na u-u\cdot\na w \rt)\cdot w\varphi^2_1dx.\nn
\eali
\ee
Using H\"older inequality, Young inequality and integration by parts, one can derive the following estimate
\bes
\bali
\int |\na w|^2\varphi^2_1dx\leq& C\lt(1+\|u\|_{L^\i(\mathfrak{C}_0)}\rt)^2 \int_{\mathfrak{C}_0} | w|^2dx\\
                \leq & C\lt(1+\|u\|_{L^\i(\mathfrak{C}_0)}\rt)^2 \|\nabla u\|^2_{L^q({\mathfrak{C}_0})},
\eali
\ees
which implies that
\bes
\|\na w\|_{L^2(\mathfrak{C}_1)}\leq \mathcal{A}\lt(\|u\|_{L^\i(\mathfrak{C}_0)},\|\nabla u\|_{L^q({\mathfrak{C}_0})} \rt).
\ees
Using Biot-Savart law and the incompressible condition, we have $-\Dl u=\na\times w$, then using the standard interior elliptic estimates, we can get
\be\label{a2ndd}
\|\na^2 u\|_{L^2(\mathfrak{C}_2)} \leq C\|\na w\|_{L^2(\mathfrak{C}_1)}\leq \mathcal{A}\lt(\|u\|_{L^\i(\mathfrak{C}_0)},\|\nabla u\|_{L^q({\mathfrak{C}_0})} \rt).
\ee
Now applying $\na$ to \eqref{erev} and then testing the resulting equation by $\na w \varphi^2_3$, we can get
\be
\bali
\int -\Dl \na w \cdot \na w\varphi^2_3dx=\int \na \lt(w\cdot\na u-u\cdot\na w\rt) \cdot \na w\varphi^2_3dx.\nn
\eali
\ee
Integration by parts, H\"{o}lder inequality and Sobolev embedding imply that
\be
\bali
&\int |\na^2 w|^2\varphi^2_3dx\\
\leq&  C\lt(1+\|u\|^2_{L^\i(\mathfrak{C}_2)}\rt) \int_{\mathfrak{C}_2}\lt( | \na w|^2+|w|^2\rt)dx+ \int _{\mathfrak{C}_2} |w\varphi^{1/2}_3|^2|\na u \varphi^{1/2}_3|^2dx\\
\leq & C\lt(1+\|u\|^2_{L^\i(\mathfrak{C}_2)}\rt) \lt(\| \na w\|^2_{L^2(\mathfrak{C}_2)}+\|w\|^2_{L^2(\mathfrak{C}_2)}\rt)+ \|w\varphi^{1/2}_3\|^2_{L^6(\mathfrak{C}_2)} \|\na u\varphi^{1/2}_3\|^2_{L^6(\mathfrak{C}_2)}\\
\leq & C\lt(1+\|u\|^2_{L^\i(\mathfrak{C}_2)}\rt) \lt(\| \na w\|^2_{L^2(\mathfrak{C}_2)}+\|w\|^2_{L^2(\mathfrak{C}_2)}\rt)\\
&+ \|\na(w\varphi^{1/2}_3)\|^2_{L^2(\mathfrak{C}_2)} \|\na(\na u\varphi^{1/2}_3)\|^2_{L^2(\mathfrak{C}_2)},\nn
\eali
\ee
which, by using \eqref{a2ndd}, implies that
\bes
\|\na^2 w\|_{L^2(\mathfrak{C}_3)} \leq \mathcal{A}\lt(\|u\|_{L^\i({\mathfrak{C}_0})},\|\nabla u\|_{L^q({\mathfrak{C}_0})} \rt).
\ees
Also the standard interior elliptic estimates indicates
\be\label{AL3}
\|\na^3 u\|_{L^2(\mathfrak{C}_4)} \leq C\|\na^2 w\|_{L^2(\mathfrak{C}_3)}\leq  \mathcal{A}\lt(\|u\|_{L^\i({\mathfrak{C}_0})},\|\nabla u\|_{L^q({\mathfrak{C}_0})} \rt).
\ee
Then applying $\na^2$ to \eqref{erev} and repeating the above procedure, similarly one can get
\be\label{AL4}
\|\na^4 u\|_{L^2(\mathfrak{C}_6)}\leq  \mathcal{A}\lt(\|u\|_{L^\i({\mathfrak{C}_0})},\|\nabla u\|_{L^q({\mathfrak{C}_0})} \rt).
\ee
Thus \eqref{eredl2} is proven by combining \eqref{a2ndd}, \eqref{AL3} and \eqref{AL4}.
\end{proof}

As a direct conclusion of Lemma \ref{A1}, one see that when $q>2$:
\be\label{erea1}
 \left\|\na w^{q/2}\right\|_{L^3(\mathcal{C}^1)}  \leq \mathcal{A}\lt(\|u\|_{L^\i({\mathcal{C}^0})},\|\nabla u\|_{L^q({\mathcal{C}^0})} \rt),
\ee
since it is easy to see that
 \be
\bali
 \left\|\na w^{q/2}\right\|_{L^3(\mathcal{C}^1)}&\ls \| w\|^{q/2-1}_{L^\i(\mathcal{C}^1)}\|\na w\|_{L^3(\mathcal{C}^1)}.\nn
\eali
\ee
Using the Sobolev embedding, we see that
\be
\| w\|_{L^\i(\mathcal{C}^1)}\ls \| (w, \na^2 w)\|_{L^2(\mathcal{C}^1)},\q \| \na w\|_{L^3(\mathcal{C}^1)}\ls \| (\na w, \na^2 w)\|_{L^2(\mathcal{C}^1)}.\nn
\ee
Then using \eqref{eredl2}, we see that \eqref{erea1} holds.
\section{Decay of $\boldsymbol{\na w^\th}$ in $\boldsymbol{L^2}$ framework}\label{appb}

 Let $\phi(x)$ be a cut-off function which satisfies
\[\text{supp}\,\phi\subset{\mathcal{C}^1}\quad\text{ and }\quad\phi(x)\equiv 1,\quad\forall\ x\in\mathcal{C}^2.\]
 Applying $\na$ to $\eqref{VEQ}_2$ and testing the resulting equation with $\na w^\th \phi^2$, it follows that
 \[
 \int_{\mathcal{C}^1}\Dl \na w^\th\na w^\th\phi^2dx=\int_{\mathcal{C}^1}\nabla\left[(b\cdot\na)w^\th+\f{1}{r^2}w^\th-\f{u^r}{r}w^\th-\f{1}{r}\p_z(u^\th)^2\right]\na w^\th \phi^2dx.
 \]
 Using integration by parts, H\"older inequality and Young inequality, noting that $r\approx 1$ in the domain $\mathcal{C}^1$, one may derive
 \[
 \begin{split}
 \int_{\mathcal{C}^1}|\na^2(w^\th\phi)|^2dx\leq&\f{1}{2}\int_{\mathcal{C}^1}|\na^2(w^\th\phi)|^2dx\\
 &+C\lt(1+\|u\|^2_{L^\i(\mathcal{C}^1)}\rt)\lt( \|\na w^\th\|^2_{L^2(\mathcal{C}^1)}+ \|\na u\|^2_{L^2(\mathcal{C}^1)}\rt),
 \end{split}
 \]
which indicates
 \be
 \|\na^2 w^\th\|^2_{L^2(\mathcal{C}^2)}\les \lt(1+\|u\|^2_{L^\i(\mathcal{C}^1)}\rt)\lt( \|\na w^\th\|^2_{L^2(\mathcal{C}^1)}+ \|\na u\|^2_{L^2(\mathcal{C}^1)}\rt).\nn
 \ee
 Using \eqref{wth1} with $q=2$, we have
 \be
 \|\na^2 w^\th\|^2_{L^2(\mathcal{C}^2)}\les \lt(1+\|u\|^3_{L^\i(\mathcal{C}^0)}\rt)\|\na u\|^2_{L^2(\mathcal{C}^0)}.\nn
 \ee
Thus the related 2-dimensional estimate follows:
\be\label{B3}
\|\bar{\na}^2 w^\th\|^2_{L^2(\mathcal{E}^2)}\les\lt(1+\|u\|^3_{L^\i(\mathcal{E}^0)}\rt) \|\bar{\na} u\|^2_{L^2(\mathcal{E}^0)}.
\ee
Noting that $\|\na^2w^\th\|_{L^3(\mathcal{E}^2)}\les\|(\nabla^3u,\na^4 u)\|_{L^2(\mathcal{E}^2)}$, applying \eqref{BGI2} to $\bar{\na} w^\th$ and using Lemma \ref{A1}, we have
\be\label{B4}
\|\bar{\na} w^\th\|_{L^\infty(\mathcal{E}^2)}\les\lt(1+\|\bar{\na} w^\th\|_{H^1(\mathcal{E}^2)}\rt)\log^{1/2}\lt(e+\mathcal{A}\lt(\|u\|_{L^\i({\mathcal{E}^0})},\|\bar{\nabla} u\|_{L^q({\mathcal{E}^0})}\rt)\rt).
\ee
Inserting \eqref{B3} \eqref{wth1} (with $q=2$) to the right hand side of \eqref{B4} and going back to the 3-dimensional domain $\mathcal{C}^2$, it follows that
\begin{small}
\be\label{B5}
\|\na w^\th\|_{L^\infty(\mathcal{C}^2)}\les\lt(1+\|u\|_{L^\infty(\mathcal{C}^0)}^{3/2}\rt)\|\nabla u\|_{L^2(\mathcal{C}^0)}\log^{1/2}\lt(e+\mathcal{A}\lt(\|u\|_{L^\i({\mathcal{C}^0})},\|\nabla u\|_{L^q({\mathcal{C}^0})}\rt)\rt).\nn
\ee
\end{small}
Now we take back the ``$\sim$" to the scaled solution, which is
\begin{small}
\be\label{B6}
\|\t{\na} \tilde{w}^\th\|_{L^\infty(\mathcal{C}^2)}\les\lt(1+\|\tilde{u}\|_{L^\infty(\mathcal{C}^0)}^{3/2}\rt)\|\nabla \tilde{u}\|_{L^2(\mathcal{C}^0)}\log^{1/2}\lt(e+\mathcal{A}\lt(\|\tilde{u}\|_{L^\i({\mathcal{C}^0})},\|\t{\nabla} \tilde{u}\|_{L^q({\mathcal{C}^0})}\rt)\rt).\nn
\ee
\end{small}
If we scale back to the domains with  ``$\lambda-$size" for $\lambda>>1$ and use H\"older inequality, then we have

\[
\begin{split}
\la^3\left\|\na{w}^\th \right\|_{L^\infty(\mathcal{C}^2_\lambda)}\lesssim&\left(1+\lambda^{3/2}\left\|u\right\|^{3/2}_{L^\infty(\mathcal{C}^0_\lambda)}\right)
\cdot\lambda^{2-3/q}\left(\int_{\mathcal{C}^0_\lambda}\left|\na u\right|^qdx\right)^{1/q}\times\\
&\log^{1/2}\lt(e+\la^M\mathcal{A}\lt(\|{u}\|_{L^\i({\mathcal{C}^0_\la})},\|{\nabla} {u}\|_{L^q({\mathcal{C}^0_\la})}\rt)\rt),
\end{split}
\]
where $M>0$ is the scaling power of $\mathcal{A}$ whose exact value is not important here since it appears inside a ``$\log$". Thus we derive
\be\label{reva1}
\begin{split}
\left\|\na{w}^\th \right\|_{L^\infty(\mathcal{C}^2_\lambda)}\lesssim&\left(1+\lambda^{3/2}\left\|u\right\|^{3/2}_{L^\infty(\mathcal{C}^0_\lambda)}\right)
\cdot\lambda^{-1-3/q}\log^{1/2}\la.
\end{split}
\ee

\end{document}